\documentclass[11pt]{article}
\usepackage{amsmath}
\usepackage{amssymb}
\usepackage{graphicx}
\usepackage{caption2}
\usepackage{psfrag}
\usepackage{amscd}
\usepackage{amsfonts}
\usepackage{bm}                                                
\usepackage{float}
\usepackage{latexsym}
\usepackage{lscape}
\usepackage{mathrsfs}
\usepackage{verbatim}
\usepackage{geometry}
\newcommand{\ds}{\displaystyle}
\newtheorem{theorem}{Theorem}[section]
\newtheorem{lemma}{Lemma}[section]

\newtheorem{corollary}{Corollary}[section]

\newcommand{\mcK}{\ensuremath{\mathcal{K}}}

\newcommand{\mcN}{\ensuremath{\mathcal{N}}}
\newcommand{\mcP}{\ensuremath{\mathcal{P}}}

\def\qed{\hbox{\vrule width 6pt height 6pt depth 0pt}}

\title{Spectral approximation of a variable coefficient fractional diffusion equation in one space dimension} 
\author{
	Xiangcheng Zheng \thanks{Department of Mathematics, University of South Carolina, Columbia,
		South Carolina 29208, USA. email: {\tt xz3@math.sc.edu \& hwang@math.sc.edu}.} 
	\and
	V.~J.~Ervin\thanks{Department of Mathematical Sciences,
	  Clemson University, Clemson, South Carolina 29634-0975, USA.
	  email: {\tt vjervin@clemson.edu}. }
	\and 
	 Hong Wang $^*$ } 

\date{\today}

\begin{document}
\maketitle

\begin{abstract}
In this article we consider the approximation of a variable coefficient (two-sided) fractional diffusion equation (FDE),
having unknown $u$. By introducing an intermediate unknown, $q$, the variable coefficient FDE is rewritten as
a lower order, constant coefficient FDE. A spectral approximation scheme, using Jacobi polynomials, is presented
for the approximation of $q$, $q_{N}$. The approximate solution to $u$, $u_{N}$, is obtained by post processing 
$q_{N}$. An a priori error analysis is given for $(q \, - \, q_{N})$ and $(u \, - \, u_{N})$. Two numerical experiments
are presented whose results demonstrate the sharpness of the derived error estimates.
\end{abstract}

\textbf{Key words}.  Fractional diffusion equation, Jacobi polynomials, spectral method

\textbf{AMS Mathematics subject classifications}. 65N30, 35B65, 41A10, 33C45

\setcounter{equation}{0}
\setcounter{figure}{0}
\setcounter{table}{0}
\setcounter{theorem}{0}
\setcounter{lemma}{0}
\setcounter{corollary}{0}
\section{Introduction}
In recent years the numerical approximation of fractional differential equations (FDEs) has received increased attention
as their incorporation into models, to address phenomena not well captured using usual differential equations, has
increased. Examples of applications using FDEs include contaminant transport in ground water flow \cite{BenWhe},
viscoelasticity \cite{mai971}, image processing \cite{bua101, gat151},
turbulent flow \cite{mai971, shl871}, and chaotic dynamics \cite{zas931}.
Approximation schemes including 
finite difference methods \cite{cui091, liu041, mee041, tad071, wan121},
finite element methods \cite{ErvRoo, jin151, liu111, WanYan}, 
discontinuous Galerkin methods \cite{xu141}, mixed methods \cite{che161, li171},
spectral methods \cite{che162, ErvHeuRoo, li121, MaoShe, mao161, MaoShe, MaoGeo, WanZha, zay151}, 
enriched subspace methods \cite{jin161}
have all been applied to FDEs.

Our interest in this paper is on the numerical approximation of the two-sided variable-coefficient FDE of order $1 < \alpha < 2$
\begin{align}
\mathcal{K}_r^{\alpha}u(x) \ := \
 -D ( \, (r \,  {}_0I_x^{2-\alpha} \ +  (1-r) \,  {}_xI_1^{2-\alpha}) \, K(x) \, Du(x))\,  &=f(x),~~x \in (0 , 1) \, ,  \label{FDE:e1a}  \\
u(0) = u(1) &= 0 \, , \label{FDE:e1b}
\end{align}
where $K(x)$ is the diffusivity coefficient with $0<K_{min}\leq K(x)\leq K_{max}$, $0\leq r\leq 1$  and $f(x)$ the source or sink term. 
The left and right fractional integrals of order $0 < \sigma < 1$ are defined as \cite{Pod,SamKil}
\begin{equation*}
{}_0I_x^{\sigma}w(x) := \frac{1}{\Gamma(\sigma)} \int_0^x \frac{w(s)}{(x-s)^{1-\sigma}} \, ds, \qquad {}_xI^{\sigma}_1w(x) := \frac{1}{\Gamma(\sigma)} \int_x^1 \frac{w(s)}{(s-x)^{1-\sigma}} \, ds \, ,
\end{equation*}
where $\Gamma(\cdot)$ is the Gamma function. Equation (\ref{FDE:e1a}) was derived by incorporating a nonlocal Fick's law with variable diffusivity coefficient $K(x)$ into a conventional local mass conservation law \cite{CasCar,ErvRoo,ZhaBen}.

In \cite{ErvRoo} the Galerkin weak formulation for (\ref{FDE:e1a}) and (\ref{FDE:e1b}) was presented and studied for $K(x)$ a constant. It was shown in \cite{WanYan} that the bilinear form of the Galerkin weak formulation may lose its coercivity for a variable-coefficient $K(x)$, and so its Galerkin finite element approximation might diverge \cite{WanYanZhu17}. A Petrov-Galerkin weak formulation was proved to be wellposed on $H^{\alpha-1}_0 \times H^1_0$ for $3/2 < \alpha < 2$ for a one-sided version of (\ref{FDE:e1a}) and (\ref{FDE:e1b}) \cite{WanYan}. A Petrov-Galerkin finite element method was developed and analyzed subsequently for the one-sided version of (\ref{FDE:e1a}) and (\ref{FDE:e1b}) \cite{WanYanZhu15}. In \cite{li171}, with the introduction of an auxiliary variable, a mixed method approximation scheme for problem (\ref{FDE:e1a}) and (\ref{FDE:e1b}) was studied and error estimates derived. 
In \cite{MaoShe}, a spectral Galerkin method for the two-sided steady-state FDE with variable coefficient was analyzed, in which the outside and inside fractional derivatives are chosen carefully so that the corresponding Galerkin weak formulation are self-adjoint and coercive. Optimal error estimates were also derived under suitable smoothness assumption on the solution. 

It was shown in \cite{WanYanZhu} that for one-dimensional FDEs that smoothness of the coefficients and the right-hand side function is not sufficient to guarantee the smoothness of the solution, especially at the endpoints of the interval, which is different from the case of the classical second order diffusion equation. Hence, seeking proper regularity solution spaces for FDEs becomes a key issue in the study of FDEs. Jin et al. \cite{jin151} conducted a thorough analysis of the regularity issue in the context of a one-sided constant-coefficient FDE by fully utilizing the explicit solution expression. An indirect Legendre spectral Galerkin method \cite{WanZha} and a finite element method \cite{WanYanZhu17} were developed for the one-sided FDE with variable coefficient, in which the solution to the FDE is expressed as a fractional derivative of the solution to a second-order differential equations. Consequently, high-order convergence rates of numerical approximations were proved using only regularity assumptions on the  coefficients and right-hand side, but not on the true solution (which is not smooth in fact). However, many aforementioned works for one-sided FDEs do not apply for two-sided FDEs.

Mao et al. \cite{mao161} analyzed the solution structure to the constant coefficient version of (\ref{FDE:e1a}) and (\ref{FDE:e1b}) with $r=1/2$ in terms of spectral polynomials and developed corresponding spectral methods. The solution structure to the constant coefficient version of (\ref{FDE:e1a}) and (\ref{FDE:e1b}) with general $0 \le r \le 1$ was resolved completely in \cite{ErvHeuRoo}, the spectral method utilizing the weighted Jacobi polynomial was studied and a priori error estimates derived. The two-sided FDE with constant coefficient and Riemann-Liouville fractional derivative was investigated in \cite{MaoGeo}, by employing a Petrov-Galerkin projection in a properly weighted Sobolev space using two-sided Jacobi polyfracnomials as test and trial functions. Spectral methods enjoy many excellent mathematical properties that make them particularly suited for FDEs: (i) They present a clean analytical expression of the true solution to FDEs, which have been fully explored in \cite{ErvHeuRoo, mao161} in analyzing the structure and regularity of the true solutions; (ii) Fractional differentiation of many spectral polynomials can be carried out analytically \cite{WanZha}, in contrast to finite element methods in which they have to be calculated numerically that are sometimes a headache \cite{WanYanZhu17}; (iii) As FDEs are nonlocal operators the appealing property of a sparse coefficient matrix, which arises for a finite element, finite difference, or finite volume approximation of a usual differential equation, is lost. In contrast, the stiffness matrices of spectral methods are often diagonal (at least for constant coefficient FDEs). Because of this, and also their convergence properties, spectral methods are appealing for the approximation of FDEs.

The goal of this paper is to extend the application of the spectral method in \cite{ErvHeuRoo} to the two-sided variable-coefficient FDE (\ref{FDE:e1a}) and (\ref{FDE:e1b}) whose solution may have endpoint singularities. By introducing an intermediate variable, we rewrite the variable coefficient model as a constant coefficient FDE. Then, utilizing Jacobi polynomials which incorporate the possible singularity of solution at endpoints, we apply the spectral method to construct a series approximation to the solution. 

This paper is orgainzed as follows. In Section \ref{secpre} we present the formulation to be used, introduce notation used through the paper, and give some key lemmas used in the analysis. The spectral approximation method is formulated and a detailed analysis of its convergence is given in Sections \ref{secspax}, \ref{secregq} and \ref{ssecEc3}. Two numerical experiments are presented in Section \ref{secNum} whose results demonstrate the sharpness of the derived error estimates.

%%%%%%%%%%%%
%%%%%%%%%%%%
%%%%%%%%%%%%

\setcounter{equation}{0}
\setcounter{figure}{0}
\setcounter{table}{0}
\setcounter{theorem}{0}
\setcounter{lemma}{0}
\setcounter{corollary}{0}
\section{Problem formulation and preliminaries}\label{secpre} 

Let $\tilde{q}(x)=-K(x)Du(x)$. Using the homogeneous Dirichlet boundary condition at $x=0$ yields
\begin{equation}\label{FDE:e2}
u(x) = -\int_0^x\frac{\tilde{q}(s)}{K(s)}ds.
\end{equation}
Enforcing the homogeneous Dirichlet boundary condition at $x=1$ we obtain 
\begin{equation}\label{FDE:e3}
\int_0^1\frac{\tilde{q}(s)}{K(s)}ds=0.
\end{equation}

Thus, with \eqref{FDE:e2}, problem (\ref{FDE:e1a}), (\ref{FDE:e1b}) can be recast as the following system
\begin{align}
\mathcal{N}_r^{\alpha}\tilde{q}(x) \, := \, 
D \, \big( r~{}_0I_x^{2-\alpha}+(1-r)~{}_xI_1^{2-\alpha}  \big) \tilde{q}(x) &= \ f(x),~~x\in (0,1) \, , \label{FDE:e4a} \\
\mathrm{with} \quad  \int_0^1\frac{\tilde{q}(s)}{K(s)}ds &= \ 0.  \label{FDE:e4b}
\end{align}

%%Important in our analysis below is the operator $\mcL^{r}_{\alpha}$ defined by
%%\begin{equation}
%%\mcL_{r}^{\alpha}  q(x) \, := \, 
%% - D^{2} \, \big( r~{}_0I_x^{2-\alpha}+(1-r)~{}_xI_1^{2-\alpha} \big) \, q(x) \, .
%% \label{defmcL}
%%\end{equation} 

Jacobi polynomial play a key role in the approximation schemes. We briefly review their definition and
properties central to the method \cite{abr641, sze751}. 

\textbf{Usual Jacobi Polynomials, $P_{n}^{(\alpha , \beta)}(x)$, on $(-1 \, , \, 1)$}.   \\    
\underline{Definition}: $ P_{n}^{(\alpha , \beta)}(x) \ := \ 
\sum_{m = 0}^{n} \, p_{n , m} \, (x - 1)^{(n - m)} (x + 1)^{m}$, where
\begin{equation}
       p_{n , m} \ := \ \frac{1}{2^{n}} \, \left( \begin{array}{c}
                                                              n + \alpha \\
                                                              m  \end{array} \right) \,
                                                    \left( \begin{array}{c}
                                                              n + \beta \\
                                                              n - m  \end{array} \right) \, .
\label{spm21}
\end{equation}
\underline{Orthogonality}:    
\begin{align}
 & \int_{-1}^{1} (1 - x)^{\alpha} (1 + x)^{\beta} \, P_{j}^{(\alpha , \beta)}(x) \, P_{k}^{(\alpha , \beta)}(x)  \, dx 
 \ = \
   \left\{ \begin{array}{ll} 
   0 , & k \ne j  \\
   |\| P_{j}^{(\alpha , \beta)} |\|^{2}
   \, , & k = j  
    \end{array} \right.  \, ,  \nonumber \\
& \quad \quad \mbox{where } \  \ |\| P_{j}^{(\alpha , \beta)} |\| \ = \
 \bigg ( \frac{2^{(\alpha + \beta + 1)}}{(2j \, + \, \alpha \, + \, \beta \, + 1)} 
   \frac{\Gamma(j + \alpha + 1) \, \Gamma(j + \beta + 1)}{\Gamma(j + 1) \, \Gamma(j + \alpha + \beta + 1)}
   \bigg )^{1/2} \, .
  \label{spm22}
\end{align}                                                    

In order to transform the domain of the family of Jacobi polynomials to $[0 , 1]$, let $x \rightarrow 2t - 1$ and 
introduce $G_{n}^{\alpha , \beta}(t) \, = \, P_{n}^{\alpha , \beta}( x(t) )$. From \eqref{spm22},
\begin{align}
 \int_{-1}^{1} (1 - x)^{\alpha} (1 + x)^{\beta} \, P_{j}^{(\alpha , \beta)}(x) \, P_{k}^{(\alpha , \beta)}(x)  \, dx 
 &= \
 \int_{t = 0}^{1} 2^{\alpha} \, (1 - t)^{\alpha} \, 2^{\beta} \, t^{\beta} \, P_{j}^{(\alpha , \beta)}(2t - 1) \, P_{k}^{(\alpha , \beta)}(2t - 1)  \, 2 \,  dt
 \nonumber \\
  &= \
2^{\alpha + \beta + 1} \int_{t = 0}^{1}   (1 - t)^{\alpha}  \, t^{\beta} \, G_{j}^{(\alpha , \beta)}(t) \, G_{k}^{(\alpha , \beta)}(t)  \,  dt
\nonumber \\
&= \
   \left\{ \begin{array}{ll} 
   0 , & k \ne j \, , \\
  2^{\alpha + \beta + 1} \, |\| G_{j}^{(\alpha , \beta)} |\|^{2}
   \, , & k = j  
   \, . \end{array} \right.    \nonumber \\
 \quad  \mbox{where } \  \  \big | \big \| G_{j}^{(\alpha , \beta)} \big | \big \| = \
 \bigg(& \frac{1}{(2j \, + \, \alpha \, + \, \beta \, + 1)} 
   \frac{\Gamma(j + \alpha + 1) \, \Gamma(j + \beta + 1)}{\Gamma(j + 1) \, \Gamma(j + \alpha + \beta + 1)}
   \bigg)^{1/2}  .  \label{spm22g} 
\end{align}                                                    

\begin{equation}
\mbox{Note that } \quad  |\| G_{j}^{(\alpha , \beta)} |\| \ = \ |\| G_{j}^{(\beta , \alpha)} |\| \, .
\label{nmeqG}
\end{equation}

 From \cite[equation (2.19)]{mao161} we have that
\begin{equation}
   \frac{d^{k}}{dx^{k}} P_{n}^{(\alpha , \beta)}(x) \ = \ 
   \frac{\Gamma(n + k + \alpha + \beta + 1)}{2^{k} \, \Gamma(n + \alpha + \beta + 1)} P_{n - k}^{(\alpha + k \, , \, \beta + k)}(x) \, .
   \label{derP}
\end{equation}   
Hence,
\begin{align}
\frac{d^{k}}{dt^{k}} G_{n}^{(\alpha , \beta)}(t) 
  &= \ \frac{\Gamma(n + k + \alpha + \beta + 1)}{  \Gamma(n + \alpha + \beta + 1)} 
  G_{n - k}^{(\alpha + k \, , \, \beta + k)}(t)  \, .  \label{eqC4}
\end{align}   

Also, from \cite[equation (2.15)]{mao161}, 
\begin{equation}
\begin{array}{l}
\ds\frac{d^{k}}{dx^{k}} \left\{ (1 - x)^{\alpha + k} \, (1 + x)^{\beta + k} \, P_{n - k}^{(\alpha + k \, , \, \beta + k)}(x) \right\}
\ \\[0.1in]
\ds\quad\quad\quad\quad\quad\quad= \ 
\frac{(-1)^{k} \, 2^{k} \, n!}{(n - k)!} \, (1 - x)^{\alpha} \, (1 + x)^{\beta} \, P_{n}^{(\alpha \, , \, \beta)}(x) \, , \
n \ge k \ge 0 \, ,
\label{eqB0}
\end{array}
\end{equation}
from which it follows that
\begin{equation}
 \frac{d^{k}}{dt^{k}} \left\{  \ (1 \, - \, t)^{\alpha + k} \,  t^{\beta + k} \,
 G_{n - k}^{(\alpha + k \, , \, \beta + k)}(t) \right\} 
 \ = \ 
 \frac{(-1)^{k} \,  n!}{(n - k)!} \,  (1 \, - \, t)^{\alpha} \,  t^{\beta} \,
 G_{n}^{(\alpha \, , \, \beta)}(t) \, . 
  \label{eqC2}
 \end{equation}
 
For compactness of notation we introduce
\begin{equation}
 \rho^{(\alpha , \beta)} \, = \, \rho^{(\alpha , \beta)}(x) \, := \, (1 - x)^{\alpha} \, x^{\beta} \, .
 \label{defrho}
\end{equation} 

 We use $y_{n} \sim n^{p}$ to denote that there exists constants $c$ and $C \ge 0$ such that, as 
 $n \rightarrow \infty$, $c \, n^{p} \le | y_{n} | \le C \, n^{p}$.

\textbf{The weighted $L^{2}(0,1)$ spaces, $L_{\omega}^{2}(0,1)$}. \\
The weighted $L^{2}(0,1)$ spaces are appropriate for analyzing the convergence of the spectral type methods 
presented below. For $\omega(x) > 0, \ x \in (0 , 1)$, let 
\[
L_{\omega}^{2}(0 , 1) \, := \, \bigg  \{ f(x) \, : \, \int_{0}^{1} \omega(x) \, f(x)^{2} \, dx \ < \ \infty  \bigg \} \, .
\]
Associated with $L_{\omega}^{2}(0 , 1)$ is the inner product, $\langle \cdot , \cdot \rangle_{\omega}$, and
norm, $\| \cdot \|_{\omega}$, defined by
\begin{align*}
\langle f \,  , \, g \rangle_{\omega} &:= \ \int_{0}^{1} \omega(x) \, f(x) \, g(x) \, dx \, , \quad \mbox{and} \\
 \| f \|_{\omega} &:= \ \left( \langle f \,  , \, f \rangle_{\omega} \right)^{1/2} \, .
\end{align*}

For $1 < \alpha < 2$ and $0 \le r \le 1$ given, let  $\beta$ satisfying $\alpha-1 \le \beta, \alpha-\beta \le 1$ be determined by
\begin{equation}
  r \ = \ \frac{\sin( \pi \, \beta)}{\sin( \pi ( \alpha - \beta)) \, + \,  \sin( \pi \, \beta)} \, . \label{defbeta} 
\end{equation}

The following two lemmas are useful in discussing the approximation scheme presented in Section~\ref{secspax}.

\begin{lemma}  \cite{ErvHeuRoo} \label{ker1r} For $\beta$ determined by \eqref{defbeta},
we have that \\
\[
ker(\mcN_{r}^{\alpha}) \ = \ span \big \{ k(x) \ := \  (1 - x)^{\alpha - \beta - 1} \, x^{\beta - 1} \big \} \, .
\]
where $\mcN_r^\alpha$ is defined in (\ref{FDE:e4a}).
%
%%Additionally, [Lue Ling Jia] (as $G_{0}^{(\delta , \gamma)}(x) \, = \, 1$)
%%\begin{align}
%% \mcN_{r}^{\alpha} (x \, k(x)) \ = \  (1 - r) \, \Gamma(\alpha) \frac{\sin(\pi \alpha)}{\sin(\pi q)}
%% &= \, \lambda_{-1} \, G_{0}^{(\delta , \gamma)}(x) \, , \quad \mbox{where } 
%%   \lambda_{-1} \, := \,  (1 - r) \, \Gamma(\alpha) \frac{\sin(\pi \alpha)}{\sin(\pi q)} \, ,
%%  \label{cherry1}  \nonumber \\
%%%
%%\mbox{and } \ \   \mcN_{r}^{\alpha} ((1 - x) \, k(x)) &= \, - \lambda_{-1} \, G_{0}^{(\alpha , \beta)}(x) \, . \nonumber
%%\end{align}  
%%\end{lemma}
%%\mbox{ } \hfill \qed   
%%
%%
%%
Additionally, \cite{che181} (as $G_{0}^{(\delta , \gamma)}(x) \, = \, 1$)
\begin{align}
  \mcN_{r}^{\alpha} (x \, k(x)) &= \  - (1 - r) \, \Gamma(\alpha) \frac{\sin(\pi \alpha)}{\sin(\pi  (\alpha - \beta))}  
 \ = \ \lambda_{-1} \, G_{0}^{(\delta , \gamma)}(x) \, ,   \nonumber \\
 \mbox{and } \ \    \mcN_{r}^{\alpha} ((1 - x) \, k(x)) &= \, - \lambda_{-1} \, G_{0}^{(\delta , \gamma)}(x) \, ,   \nonumber \\
 \quad \mbox{where } 
   \lambda_{-1} \, &:= \,  - (1 - r) \, \Gamma(\alpha) \frac{\sin(\pi \alpha)}{\sin(\pi (\alpha - \beta))} \, .
  \nonumber 
\end{align}  
\end{lemma}
\mbox{ } \hfill \qed

In the representation of $q_{N}(x)$ (see \eqref{defuqNs}), the approximation of $q(x)$, we need to include either
$x \, k(x)$ or $(1 - x) \, k(x)$. If $\beta - 1 \,  < \, \alpha - \beta$, 
i.e., $r > 1/2$, then $x \, k(x)$ is a more regular function on $(0 , 1)$ than
$(1 - x) \, k(x)$. However, if  $\beta - 1 \,  > \, \alpha - \beta$, i.e., 
$r < 1/2$, then $(1 - x) \, k(x)$ is a more regular function on $(0 , 1)$ than
$x \, k(x)$. 

\begin{lemma} \cite{ErvHeuRoo} \label{lmaegp}
Let $\beta$ be determined by \eqref{defbeta}.
Then, for $n \, = \, 0, 1, 2, \ldots$
\begin{align}
 \mcN_{r}^{\alpha} (1 - x)^{\alpha - \beta} x^{\beta} \, G_{n}^{(\alpha - \beta \, , \, \beta)}(x) 
 &= \ \lambda_{n}  \, G_{n + 1}^{(\beta - 1 \, , \, \alpha - \beta - 1)}(x) \, ,  \label{NactG}  \\
 \mbox{where }  \quad
   \lambda_{n} &= \ \frac{\sin(\pi \alpha)}{\sin(\pi (\alpha - \beta)) \ + \ \sin(\pi \beta)} \, \frac{\Gamma(n + \alpha)}{n !} \, .
 \label{defln}
\end{align}
\end{lemma}
\mbox{ } \hfill \qed

Using Stirling's formula we have that
\begin{equation}
\lim_{n \rightarrow \infty} \, \frac{\Gamma(n + \mu)}{\Gamma(n) \, n^{\mu}}
\ = \ 1 \, , \mbox{ for } \mu \in \mathbb{R}.  
 \label{eqStrf}
\end{equation} 
\begin{equation}
 \mbox{Thus $\lambda_{n} > 0$ for all $n \, = \, 0, 1, 2, \ldots $, and 
$\lambda_{n} \sim   (n + 1)^{\alpha - 1}$.}
\label{nwbde}
\end{equation}

%%%%%
%%%%%
%%%%%
\setcounter{equation}{0}
\setcounter{figure}{0}
\setcounter{table}{0}
\setcounter{theorem}{0}
\setcounter{lemma}{0}
\setcounter{corollary}{0}
\section{Spectral type approximation to \eqref{FDE:e4a}, \eqref{FDE:e4b}.}
\label{secspax}
In this section we fix the values of $\alpha$ and $r$ as defined by the operator $\mcK_{r}^{\alpha}$ in \eqref{FDE:e1a},
and correspondingly, $\beta$ determined by \eqref{defbeta}. 
We will assume that $r \ge 1/2$. Hence we include $x \, k(x)$ (and not $(1 - x) \, k(x)$)
in the representation of $q_{N}(x)$ (see \eqref{defuqNs}). 

Useful in the analysis below is the following results.
\begin{lemma} \label{lmageq}
For $j \, = \, 0, 1, 2, \ldots$
\begin{equation}
 \frac{1}{2} \, \le \, \frac{ |\| G_{j}^{(\alpha - \beta \, , \, \beta)} |\|^{2} }{ |\| G_{j + 1}^{(\beta - 1 \, , \, \alpha - \beta - 1)} |\|^{2} } 
  \, = \, \frac{j + 1}{j + \alpha} \, \le \, 1 \, .
 \label{gge2}
\end{equation}
\end{lemma}
\textbf{Proof}: From \eqref{spm22g},
\begin{align}
\frac{ |\| G_{j}^{(\alpha - \beta \, , \, \beta)} |\|^{2} }{ |\| G_{j + 1}^{(\beta - 1 \, , \, \alpha - \beta - 1)} |\|^{2} } 
&= \ \frac{1}{2 j \, + \, \alpha \, + \, 1} \, 
\frac{\Gamma(j + \alpha - \beta + 1) \, \Gamma(j + \beta + 1)}{\Gamma(j + 1) \, \Gamma(j + \alpha + 1)} \nonumber\\
& \quad  \quad \quad \quad \quad \times\frac{2 j \, + \, \alpha \, + \, 1}{1} \, 
\frac{\Gamma(j + 2) \, \Gamma(j + \alpha)}{\Gamma(j + \beta + 1) \, \Gamma(j + \alpha - \beta + 1)}   
= \ \frac{j + 1}{j + \alpha} \le 1 \, .      \label{nrto1}
\end{align}
\mbox{ } \hfill \qed

The solution $u(x)$ to \eqref{FDE:e1a}, \eqref{FDE:e1b} is computed directly using \eqref{FDE:e2} once
$\tilde{q}(x)$ satisfying  \eqref{FDE:e4a}, \eqref{FDE:e4b} is determined. Note that as 
$ker(\mcN_{r}^{\alpha}) \ = \ span\{ k(x) \}$, then $\tilde{q}$ satisfying \eqref{FDE:e4a}
is only determined up to an additive constant multiple of $k(x)$. Hence, we rewrite
$\tilde{q}(x) \ = \ c_{-2} \, k(x) \ + \ q(x)$, where $q(x)$ satisfied 
\begin{equation}
\mcN_{r}^{\alpha} q(x) \ = \ f(x) , \ \ x \in (0,1) \, ,
\label{ddqq}
\end{equation}
and $c_{-2}$ is determined by \eqref{FDE:e4b}.

\textbf{Remark}: Note that $f(x) \in L_{\rho^{(\beta - 1 \, , \, \alpha - \beta - 1)}}^{2}(0,1)$ may be expressed as \\
$f(x) \ = \ \sum_{i = 0}^{\infty} \frac{f_{i}}{ |\| G_{i}^{(\beta - 1 \, , \alpha - \beta - 1)} |\|^{2}} \, 
G_{i}^{(\beta - 1 \, , \alpha - \beta - 1)}(x)$, 
where $f_{i}$ is given by
\begin{equation}
f_{i} \, := \, \int_{0}^{1} \,  \rho^{(\beta - 1\, , \, \alpha - \beta - 1)}(x) \, f(x) \, G_{i}^{(\beta - 1 \, , \alpha - \beta - 1)}(x)  \, dx \, .
\label{deffis}
\end{equation}

With $f_{i}$ defined in \eqref{deffis}, let 
\begin{align}
q_{N}(x)  &= \  c_{-1} x \, k(x) \ + \ \rho^{(\alpha - \beta \, , \, \beta)}(x) \, 
\sum_{i = 0}^{N-1} c_{i} \, G_{i}^{(\alpha - \beta \, , \, \beta)}(x) \, ,   \label{defuqNs} \\
\mbox{ where} \ &
c_{i} \ = \ \frac{1}{ \lambda_{i} \, |\| G_{i+1}^{(\beta - 1 \, , \alpha - \beta - 1)} |\|^{2}} f_{i+1} \, , \mbox{for } i = -1, 0, 1, \ldots N-1.
\nonumber 
\end{align}

\begin{theorem}   \label{thmuspg}
Let $f (x) \in L_{\rho^{(\beta - 1 \, , \, \alpha - \beta - 1)}}^{2}(0,1)$ and $q_{N}(x)$ be as defined in \eqref{defuqNs}. Then, \\
$\left( q(x) -  c_{-1} x \, k(x) \right) \ := \ \lim_{N \rightarrow \infty} \left( q_{N}(x) -  c_{-1} x \, k(x) \right) \in L^{2}_{\rho^{(-(\alpha - \beta) \, , - \beta)}}(0,1)$. In addition, \\
$\mcN_{r}^{\alpha} q(x) \ = \ f(x)$.
\end{theorem}
\textbf{Proof}: For $f_{N}(x) \ := \ \sum_{i = 0}^{N} 
\frac{f_{i}}{ |\| G_{i}^{(\beta - 1 \, , \alpha - \beta - 1)} |\|^{2}} \, 
G_{i}^{(\beta - 1 \, , \alpha - \beta - 1)}(x)$, we have that $f(x) \ = \ \lim_{N \rightarrow \infty} f_{N}(x)$, and 
$\left\{ f_{N}(x) \right\}_{N = 0}^{\infty}$ is a Cauchy sequence in $L^{2}_{\rho^{(\beta - 1 \, , \, \alpha - \beta - 1)}}(0,1)$. 
A straightforward calculation shows that $ \left( q_{N}(x) -  c_{-1} x \, k(x) \right) \, \in L^{2}_{\rho^{(-(\alpha - \beta) \, , - \beta)}}(0,1)$.
Then, (without loss of 
generality, assume $M > N$)
\begin{align*}
& \| \left( q_{M}(x) -  c_{-1} x \, k(x) \right) - \  \left( q_{N}(x) -  c_{-1} x \, k(x) \right) \|_{\rho^{(-(\alpha - \beta) \, , - \beta)}}^{2} \ = \ 
\| q_{M}(x) \ - \ q_{N}(x) \|_{\rho^{(-(\alpha - \beta) \, , - \beta)}}^{2}  \\
\quad &= \ \bigg( \rho^{(-(\alpha - \beta) \, , - \beta)}(x) \ \rho^{(\alpha - \beta \, , \, \beta)}(x) \, 
\sum_{j = N}^{M - 1} c_{j} \, G_{j}^{(\alpha - \beta \, , \, \beta)}(x) \ , \ 
  \rho^{(\alpha - \beta \, , \, \beta)}(x) \, \sum_{j = N}^{M - 1} c_{j} \, G_{j}^{(\alpha - \beta \, , \, \beta)}(x) \bigg)     \\
\quad &= \
\bigg(   \rho^{(\alpha - \beta \, , \, \beta)}(x) \, \sum_{j = N}^{M - 1} \frac{G_{j}^{(\alpha - \beta \, , \, \beta)}(x)}{\lambda_{j} \, |\| G^{(\beta - 1 \, , 
\, \alpha - \beta - 1)}_{j + 1}  |\|^{2} }  \, 
 f_{j + 1}\ , \ 
   \sum_{j = N}^{M - 1} \frac{G_{j}^{(\alpha - \beta \, , \, \beta)}(x)}{\lambda_{j} \, |\| G^{(\beta - 1 \, , \, \alpha - \beta - 1)}_{j + 1}  |\|^{2} }  
   \, f_{j + 1} \bigg)  \\
\quad &= \
 \sum_{j = N}^{M - 1} \frac{f_{j + 1}^{2}}{\lambda_{j}^{2} \, |\| G^{(\beta - 1 \, , \, \alpha - \beta - 1)}_{j + 1}  |\|^{4} } 
   |\| G_{j}^{(\alpha - \beta \, , \, \beta)}  |\|^{2} \le \
 \sum_{j = N}^{M - 1} \frac{f_{j + 1}^{2}}{\lambda_{j}^{2} \, |\| G^{(\beta - 1 \, , \, \alpha - \beta - 1)}_{j + 1} |\|^{2} }   
  \quad \mbox{ (using \eqref{gge2})} \\
\quad &\le \ C \,
\bigg(  \rho^{(\beta - 1 \, , \, \alpha - \beta - 1)}(x) \, \sum_{j = N + 1}^{M} \frac{G^{(\beta - 1 \, , \, \alpha - \beta - 1)}_{j}(x)}{ |\| G^{(\beta - 1 \, , \, \alpha - \beta - 1)}_{j}  |\|^{2} }  f_{j}\, 
 \ , \ 
   \sum_{j = N + 1}^{M} \frac{ G^{(\beta - 1 \, , \, \alpha - \beta - 1)}_{j}(x)}{ |\| G^{(\beta - 1 \, , \, \alpha - \beta - 1)}_{j}  |\|^{2} } f_{j} \,
    \bigg)  \\
\quad & \quad \quad \mbox{ (using $\lambda_{j}$'s are bounded away from zero) }  \\
\quad &= \ C \, \| f_{N}(x) \ - \ f_{M}(x) \|_{\rho^{(\beta - 1 \, , \, \alpha - \beta - 1)}}^{2} \, .
\end{align*}
Hence $\{ \left( q_{N}(x) -  c_{-1} x \, k(x) \right) \}_{N = 0}^{\infty}$ is a Cauchy sequence in 
$L^{2}_{\rho^{(-(\alpha - \beta) \, , - \beta)}}(0,1)$. As 
$L^{2}_{\rho^{(-(\alpha - \beta) \, , - \beta)}}(0,1)$ is complete \cite{hes071},   
$q(x) -  c_{-1} x \, k(x) \ := \ \lim_{N \rightarrow \infty} q_{N}(x) -  c_{-1} x \, k(x) \, \in \, L^{2}_{\rho^{(-(\alpha - \beta) \, , - \beta)}}(0,1)$.

Next, as $f_{N}(x) \rightarrow f(x)$ in $L^{2}_{\rho^{(\beta - 1 \, , \, \alpha - \beta - 1)}}(0 , 1)$, given $\epsilon > 0$ there exists $\tilde{N}$ 
such that for $N > \tilde{N}$, $\| f(x)  \ - \ f_{N}(x) \|_{\rho^{(\beta - 1 \, , \, \alpha - \beta - 1)}} \, < \, \epsilon$. Then, for $N > \tilde{N}$,
using Lemmas \ref{ker1r} and \ref{lmaegp}
\begin{align*}
& \| f(x) \ - \ \mcN^{\alpha}_{r} q_{N}(x) \|_{\rho^{(\beta - 1 \, , \, \alpha - \beta - 1)}} \\
\quad  &=  \ \bigg\| f(x) \ - \  
\mcN^{\alpha}_{r} \bigg(  c_{-1} x \, k(x) \ + \  \rho^{(\alpha - \beta \, , \, \beta)}(x) \, \sum_{j = 0}^{N - 1} 
 \frac{G_{j}^{(\alpha - \beta \, , \, \beta)}(x)}{\lambda_{j} \, |\| G^{(\beta - 1 \, , \, \alpha - \beta - 1)}_{j + 1}  |\|^{2} }
  f_{j + 1}\,  \bigg) \bigg\|_{\rho^{(\beta - 1 \, , \, \alpha - \beta - 1)}} \\
\quad &= \ 
\bigg\| f(x) \ - \  
\sum_{j = 0}^{N}  \frac{ G^{(\beta - 1 \, , \, \alpha - \beta - 1)}_{j}(x)}{ |\| G^{(\beta - 1 \, , \, \alpha - \beta - 1)}_{j}  |\|^{2} }
 \,f_{j}  \bigg\|_{\rho^{(\beta - 1 \, , \, \alpha - \beta - 1)}} \\
\quad &= \ 
\| f(x) \ - \  f_{N}(x)  \|_{\rho^{(\beta - 1 \, , \, \alpha - \beta - 1)}} \ < \ \epsilon \, .
\end{align*}
Hence, $f(x) \ = \ \mcN^{\alpha}_{r} q(x)$. \\
\mbox{ } \hfill \qed

For $q \, - \, q_{N}$ we have the following a priori error estimate.
\begin{theorem}
For $f(x) \in L^{2}_{\rho^{(\beta - 1 \, , \, \alpha - \beta - 1)}}(0,1)$ and $q_{N}(x)$ given by  \eqref{defuqNs}, there exists $C > 0$ such that
\begin{equation}
 \| q \, - \, q_{N} \|_{\rho^{(-(\alpha - \beta) \, , - \beta)}} \ \le \   \frac{1}{ \lambda_{N} } \, \| f \|_{\rho^{(\beta - 1 \, , \, \alpha - \beta - 1)}} \ \le \
 \ C \, (N + 1)^{- \alpha + 1}    \, \| f \|_{\rho^{(\beta - 1 \, , \, \alpha - \beta - 1)}} . 
  \label{apertg2}
\end{equation} 
\label{thmAP2}
\end{theorem} 
\textbf{Proof}: Using the definition of the $ \| \cdot \|_{\rho^{(-(\alpha - \beta) \, , - \beta)}} $ norm, 
\begin{align*}
 \| q \, - \, q_{N} \|_{\rho^{(-(\alpha - \beta) \, , - \beta)}}^{2} &=
 \int_{0}^{1} \rho^{(-(\alpha - \beta) \, , - \beta)}(x) \bigg( \rho^{(\alpha - \beta \, , \, \beta)}(x) \, \sum_{i = N}^{\infty}  
 \frac{G^{(\alpha - \beta \, , \, \beta)}_{i}(x)}{ ( \lambda_{i} \,  |\| G^{(\beta - 1 \, , \, \alpha - \beta - 1)}_{i+1} |\|^{2} )}
  \,   f_{i+1}  \bigg)^{2} \, dx  \\
%
%%&\le \ \max_{i \, \ge \, N} \left( \frac{1}{  \lambda_{i}^{2}  } \right) \ 
%%  \int_{0}^{1} \, \rho(x) \, \left( \sum_{i \, = \, N + 1}^{\infty}  \frac{ f^{*}_{i+1} }{\| G^{(\beta - 1 \, , \, \alpha - \beta - 1)}_{i+1} \|_{\varpi}^{2} } \,
%%   \mcG^{\rho}_{i}(x)   \right)^{2} \, dx  \\
%
&\le \ \max_{i \, \ge \, N} \left( \frac{1}{  \lambda_{i}^{2}  } \right)  \ 
  \sum_{i \, = \, N }^{\infty}   \frac{ f^{2}_{i + 1} }{|\| G^{(\beta - 1 \, , \, \alpha - \beta - 1)}_{i+1} |\|^{4} } \, 
  |\| G^{(\alpha - \beta \, , \, \beta)}_{i} |\|^{2} 
   \\
&\le \  \frac{1}{  \lambda_{N}^{2}  }  \ 
  \sum_{i \, = \, N}^{\infty}  \frac{ f^{2}_{i + 1} }{|\| G^{(\beta - 1 \, , \, \alpha - \beta - 1)}_{i + 1} |\|^{4} } \,
   |\| G^{(\beta - 1 \, , \, \alpha - \beta - 1)}_{i + 1} |\|^{2} 
  \quad \mbox{ (using \eqref{gge2})}
   \\
&\le \   \frac{1}{  \lambda_{N}^{2}  }   \ 
  \int_{0}^{1} \, \rho^{(\beta - 1 \, , \, \alpha - \beta - 1)}(x) \,
  \left( \sum_{i \, = \, 0}^{\infty}  \frac{ G^{(\beta - 1 \, , \, \alpha - \beta - 1)}_{i}(x) }{|\| G^{(\beta - 1 \, , \, \alpha - \beta - 1)}_{i} |\|^{2} } \, 
   f_{i}
  \right)^{2} \, dx  \\
&= \ \frac{1}{  \lambda_{N}^{2}  }   \
  \int_{0}^{1} \, \rho^{(\beta - 1 \, , \, \alpha - \beta - 1)}(x) \, f(x)^{2} \, dx  \le \  \frac{1}{  \lambda_{N}^{2}  }   \ 
  \| f \|_{\rho^{(\beta - 1 \, , \, \alpha - \beta - 1)}}^{2}  \\
&\le \ C \, (N + 1)^{-2 \, (\alpha \, - \, 1)}   \ 
  \| f \|_{\rho^{(\beta - 1 \, , \, \alpha - \beta - 1)}}^{2} \, , \quad \mbox{ using } \eqref{nwbde} .
\end{align*}
\mbox{ } \hfill \qed

\begin{corollary}
For $f(x) \in L^{2}_{\rho^{(\beta - 1 \, , \, \alpha - \beta - 1)}}(0,1)$ and $q_{N}(x)$ 
given by  \eqref{defuqNs}, there exists $C > 0$ such that
\begin{equation}
 \| q \, - \, q_{N} \|  \ \le \   \frac{1}{\lambda_{N} } \, \| f \|_{\rho^{(\beta - 1 \, , \, \alpha - \beta - 1)}} \ \le \
 \ C \, (N + 1)^{- \alpha + 1 }    \, \| f \|_{\rho^{(\beta - 1 \, , \, \alpha - \beta - 1)}} \, . \label{apertg3}
\end{equation}
\end{corollary} 
\label{corAPg1}
\textbf{Proof}: As $\rho^{(-(\alpha - \beta) \, , - \beta)}(x) \, = \, (1 - x)^{-(\alpha - \beta)} \, x^{-\beta}  \, > 1 $, for $0 < x < 1$, then
$\| u \, - \, u_{N} \| \ \le \ \| u \, - \, u_{N} \|_{\rho^{(-(\alpha - \beta) \, , - \beta)}}$. Hence 
the bound \eqref{apertg3} follows immediately from \eqref{apertg2} .
  \\
\mbox{ } \hfill \qed 

%%%%
%%%%
%%%%

\setcounter{equation}{0}
\setcounter{figure}{0}
\setcounter{table}{0}
\setcounter{theorem}{0}
\setcounter{lemma}{0}
\setcounter{corollary}{0}
%
% Herein we investigate the regularity of q_{N} - singular part
%
\section{Regularity of $D^{j}\big( (q \ - \ c_{-1} x k(x))/\rho^{(\alpha-\beta,\beta)}(x)\big)$} 
\label{secregq}
In this section we investigate the regularity of $(q \ - \ c_{-1} x k(x))/\rho^{(\alpha-\beta,\beta)}(x)$. We do this by establishing that 
$\{ D^{j} \big( (q_{N} \ - \ c_{-1} x k(x))/\rho^{(\alpha-\beta,\beta)}(x) \big) \}$ is a Cauchy sequence in an appropriately weighted $L^{2}$ function space.

Let 
\begin{align}
f_{N}(x) &= \ \sum_{i = 0}^{N} \frac{f_{i}}{ |\| G_{i}^{(\beta - 1 \, , \, \alpha - \beta - 1)}(x) |\|^2 }\,
  G_{i}^{(\beta - 1 \, , \, \alpha - \beta - 1)}(x) \, .  \nonumber \\
\mbox{Hence, using \eqref{eqC4} and } & \mbox{reindexing}    \ \   \nonumber \\
D^{j} \, f_{N}(x) &= \ \sum_{i = -1}^{N-1} \frac{f_{i + 1}}{ |\| G_{i + 1}^{(\beta - 1 \, , \, \alpha - \beta - 1)} |\|^2 } \,
\frac{\Gamma(i + j + \alpha)}{  \Gamma(i + \alpha)}
  G_{i - j + 1}^{(\beta + j - 1 \, , \, \alpha - \beta + j - 1)}(x) \, .  \label{eqC5}  
\end{align}  

Helpful in establishing the general result is the following lemma.
\begin{lemma} \label{lmahf1}
For $j \in \mathbb{N}$, there exists $C > 0$ such that
\begin{equation}
 \frac{1}{\lambda_{i}^{2}} \, \left( \frac{i + j + \alpha}{i + \alpha} \right)^{2} \, 
 \frac{ |\| G_{i - j}^{(\alpha - \beta + j \, , \, \beta + j)} |\|^{2}}{ |\| G_{i - j + 1}^{(\beta + j - 1 \, , \, \alpha -  \beta + j - 1)} |\|^{2}}
  \ \le \ C \,  i^{- 2(\alpha - 1)}  .
 \label{errt2}
\end{equation}
\end{lemma}
\textbf{Proof}: From \eqref{nmeqG} and \eqref{spm22g} ,
\begin{align}
  \frac{ |\| G_{i - j}^{(\alpha - \beta + j \, , \, \beta + j)} |\|^{2}}{ |\| G_{i - j + 1}^{(\beta + j - 1 \, , \, \alpha -  \beta + j - 1)} |\|^{2}}
 &= \ 
 \frac{ |\| G_{i - j}^{(\alpha - \beta + j \, , \, \beta + j)} |\|^{2}}{ |\| G_{i - j + 1}^{(\alpha -  \beta + j - 1 \, , \, \beta + j - 1)} |\|^{2}}
 \nonumber \\
 & = \  \frac{1}{(2 i \,  + \,  \alpha \, + \, 1)} \, 
  \frac{\Gamma( i + \alpha - \beta + 1) \, \Gamma(i + \beta +1)}{\Gamma(i - j + 1) \, \Gamma(i + j + \alpha + 1)}  
  \nonumber \\
& \quad \quad \quad 
  \cdot \, (2 i \, + \, \alpha \, + \, 1) \, 
  \frac{\Gamma( i  - j + 2) \, \Gamma(i + j + \alpha)}{\Gamma(i + \alpha - \beta  + 1) \, \Gamma(i + \beta + 1) }  
  \nonumber \\
&= \ 
  \frac{(i - j + 1)}{(i + j + \alpha)}  . \label{ppol1}
%%   \nonumber \\
%%%
%%&\sim \frac{\Gamma(i \, - \, 2 k \, + 2)}{\Gamma(i + 1)} \, \frac{\Gamma(i + \alpha + 2)}{\Gamma(i + \alpha \, - \, 2 k \, + 1)} \,
%%  \frac{\Gamma( i + \beta - k + 1)}{\Gamma(i + \beta - k + 3)} \, \ \  \mbox{ (for $i$ large)}  \,  \nonumber  \\
%%%
%%&= \frac{1}{i \, (i - 1) \ldots (i - (2k - 2))} \frac{\Gamma(i \, - \, 2 k \, + 2)}{\Gamma(i \, - \, 2 k \, + 2)} \nonumber \\
%%&  \quad \quad  \cdot  \frac{(i + \alpha + 1) (i + \alpha) \ldots (i + \alpha + 1 - 2k)}{1} 
%%  \frac{\Gamma(i + \alpha \, - \, 2 k \, + 1)}{\Gamma(i + \alpha \, - \, 2 k \, + 1)} \nonumber \\
%%& \quad \quad  \quad \quad  \cdot \frac{1}{(i + \beta - k + 2) \, (i + \beta - k + 1)} 
%%   \frac{\Gamma( i + \beta - k + 1)}{\Gamma(i + \beta - k + 1)}  \nonumber \\
%%%
%%&\sim i^{- (2k - 1)} \ i^{2k + 1} \ i^{- 2} \ = \ 1 \, .   \label{ppol1}
\end{align}

Using Stirling's formula,
\begin{align}
 \frac{1}{| \lambda_{i} |} \ = \ C \ \frac{\Gamma(i + 1)}{\Gamma(i + \alpha)}
  &\sim 
   \ \left( i + 1 \right)^{- (\alpha - 1)} \ \sim \ i^{- (\alpha - 1)}  \, . \label{uy3}   \end{align}
Combining \eqref{ppol1} and \eqref{uy3} we obtain
\begin{align*}
\frac{1}{\lambda_{i}^{2}} \, \left( \frac{i + j + \alpha}{i + \alpha} \right)^{2}  \, 
 \frac{ |\| G_{i - j}^{(\alpha - \beta + j \, , \, \beta + j)} |\|^{2}}{ |\| G_{i - j + 1}^{(\beta + j - 1 \, , \, \alpha -  \beta + j - 1)} |\|^{2}} 
 \sim 
 \left( i^{-(\alpha - 1)} \right)^{2} \,  \left( \frac{i + j + \alpha}{i + \alpha} \right)^{2}  \,
   \frac{(i - j + 1)}{(i + j + \alpha)} \  \sim \ i^{- 2(\alpha - 1)}   \, ,
\end{align*}
from which \eqref{errt2} follows.  \\
\mbox{ } \hfill \qed

We have the following theorem.
\begin{theorem} \label{thmrgqd}
For $j \in \mathbb{N}$, if $D^{j}f \in L^{2}_{\rho^{(\beta + j - 1 \, , \, \alpha - \beta + j - 1)}}(0,1)$, then
$D^{j}\big((q(x) \, - \, c_{-1} x \, k(x))/\rho^{(\alpha-\beta,\beta)}(x)\big) \in L^{2}_{\rho^{(\alpha - \beta + j \, , \, \beta + j )}}(0,1)$.
\end{theorem}
\textbf{Proof}:
From \eqref{eqC2} and \eqref{defuqNs},
\begin{align}
D^{j} \left( \frac{q_{N} \ - \ c_{-1} x \,  k(x)}{\rho^{(\alpha-\beta,\beta)}(x)} \right)
&= \ D^{j} \bigg(  \sum_{i = 0}^{N - 1} c_{i} \, G_{i}^{(\alpha - \beta \, , \, \beta)}(x) \bigg)=  \sum_{i = 0}^{N - 1 } c_{i } \, 
\frac{\Gamma(i+j+\alpha+1)}{\Gamma(i+\alpha+1)} \, G_{i-j}^{(\alpha - \beta + j \, , \, \beta + j)}(x) ,  \nonumber 
\end{align}
where $G_{k}^{(a , b)}(x) = 0$ for $k < 0$.

Then, 
\begin{align}
&
\left\| D^{j} \left( \frac{q_{M} \ - \ c_{-1} xk(x)}{\rho^{(\alpha-\beta,\beta)}(x)} \right) \ - \ D^{j} \left( \frac{q_{N} \ - \ c_{-1}x k(x)}{\rho^{(\alpha-\beta,\beta)}(x)} \right) \right\|^{2}_{\rho^{(\alpha - \beta + j \, , \, \beta + j )}}
 \nonumber  \\
& \ = \ \bigg( \rho^{(\alpha - \beta + j \, , \, \beta +j)} \, \sum_{i = N}^{M-1} c_{i} \, 
\frac{\Gamma(i+j+\alpha+1)}{\Gamma(i+\alpha+1)} \, G_{i-j}^{(\alpha - \beta +j \, , \, \beta + j)} \ , \ 
  \sum_{i = N}^{M-1} c_{i} \, 
\frac{\Gamma(i+j+\alpha+1)}{\Gamma(i+\alpha+1)} \, G_{i-j}^{(\alpha - \beta +j \, , \, \beta + j)} \bigg)    \nonumber \\
&= \ 
 \sum_{i =N}^{M-1} c_{i}^{2} \, 
\left( \frac{\Gamma(i+j+\alpha+1)}{\Gamma(i+\alpha+1)} \right)^{2} \, |\| G_{i-j}^{(\alpha - \beta + j \, , \, \beta + j)}(x) |\|^{2}  \nonumber \\
&= \ 
 \sum_{i =N}^{M-1} \frac{f_{i+1}^2}{\lambda_i^2|\|G_{i+1}^{(\beta-1,\alpha-\beta-1)}(x)|\|^4} \, 
\left( \frac{\Gamma(i+j+\alpha+1)}{\Gamma(i+\alpha+1)} \right)^{2} \, |\| G_{i-j}^{(\alpha - \beta + j \, , \, \beta + j)}(x) |\|^{2}  \nonumber \\
&\le \ 
C   \sum_{i =N}^{M-1} \frac{f_{i+1}^2}{|\|G_{i+1}^{(\beta-1,\alpha-\beta-1)}(x)|\|^4} \, i^{- 2(\alpha - 1)} \, 
\left( \frac{\Gamma(i+j+\alpha)}{\Gamma(i+\alpha)} \right)^{2} \, |\| G_{i+1-j}^{(\beta + j-1 \, , \, \alpha - \beta + j-1)}(x) |\|^{2}  \ \mbox{ (using \eqref{errt2})} \nonumber \\
&\le C \, N^{- 2(\alpha - 1)}  \, \| D^{j} f_{M}(x) \, - \, D^{j} f_{N}(x) \|_{\rho^{(\beta + j - 1 \, , \, \alpha - \beta + j - 1)}}^{2}  \ \ \ \ 
\mbox{(using \eqref{eqC5})}  \, ,  \label{thy1}  \\
&= C \,  \| D^{j} f_{M}(x) \, - \, D^{j} f_{N}(x) \|_{\rho^{(\beta + j - 1 \, , \, \alpha - \beta + j - 1)}}^{2} \, . \nonumber
\end{align}

Assuming that $D^{j}f \in L^{2}_{\rho^{(\beta + j - 1 \, , \, \alpha - \beta + j - 1)}}(0,1)$, then $\{D^{j} f_{n} \}$ is a Cauchy sequence in 
$L^{2}_{\rho^{(\beta + j - 1 \, , \, \alpha - \beta + j - 1)}}(0,1)$. Thus we can conclude that 
$D^{j} \big((q \, - \, c_{-1} x \, k(x))/\rho^{(\alpha-\beta,\beta)}(x)\big ) \in  
L^{2}_{\rho^{(\alpha - \beta + j \, , \, \beta + j )}}(0,1)$.  \\
\mbox{ } \hfill \qed
 
%%%%%%%
%%%%%%%
\subsection{Additional error estimate for $q \ - \ c_{-1}k(x)$} 
\label{secAddErr}
From Theorems \ref{thmuspg} and \ref{thmrgqd} we have that  
$q \ - \ c_{-1} x  \, k(x) \in L^{2}_{\rho^{(-(\alpha - \beta) \, , \, - \beta)}}(0,1)$,
and $D^{j}\big((q(x) \, - \, c_{-1} x \, k(x))/\rho^{(\alpha-\beta,\beta)}(x)\big) \in 
L^{2}_{\rho^{(\alpha - \beta + j \, , \, \beta + j )}}(0,1)$, $j \in \mathbb{N}$, for 
a sufficiently smooth rhs function, $f(x)$. Thus, for each successive derivative of 
$(q(x) \, - \, c_{-1} x \, k(x))/\rho^{(\alpha-\beta,\beta)}$ the power of the weight function
at the endpoints of the interval needs to be increased by one for the function to be (weighted) square integrable. This observation
leads to the following definition of weighted Sobolev spaces \cite{guo041}.
\[
  H^{r}_{\rho^{(a \, , \, b)} \,  , \, A}(0,1) \ := \
  \left\{ v \, | \, v \mbox{ is measurable and } \| v \|_{r , \, \rho^{(a \, , \, b)} , \, A} < \infty \right\} \, , \ \ r \in \mathbb{N} \, ,
\]
with associated norm and semi-norm
\[
   \| v \|_{r , \, \rho^{(a \, , \, b)} , \, A}  \ := \ \bigg( \sum_{j = 0}^{r} 
   \| D^{j} v \|_{\rho^{(a + j \, , b + j)}}^{2} \bigg)^{1/2} \, , \ \ \
    | v |_{r , \, \rho^{(a \, , \, b)} , \, A}  \ := \  \| D^{r} v \|_{\rho^{(a + r \, , b + r)}} \, .
\]

Let $\mcP_{N}$ denote the space of polynomials of degree $\le N$, and introduce the orthogonal projection
$P_{N , a , b} \, : \, L^{2}_{\rho^{(a \, , \, b)}}(0,1) \rightarrow \mcP_{N}$ defined by
\[
 \left( v \, - \,  P_{N , a , b} v \ , \ \phi \right)_{\rho^{(a \, , \, b)}} \ = \ 0 \, , \ \ \forall \phi \in \mcP_{N} \, .
\]
Then from \cite{guo041} we have the following theorem.
\begin{theorem} \cite[Theorem 2.1]{guo041} 
For any $v \in H^{r}_{\rho^{(a \, , \, b)} \,  , \, A}(0,1)$, $r \in \mathbb{N}$, and $0 \le \mu \le r$, there exists a
constant $C$, independent of $N, \, \alpha$ and $\beta$ such that
\begin{equation}
  \| v(x) \, - \,  P_{N , a , b} v(x) \|_{\mu , \, \rho^{(a \, , \, b)} , \, A} \ \le \ C \, 
  ( N \, (N + a + b) )^{\frac{\mu - r}{2}} \, | v |_{r , \, \rho^{(a \, , \, b)} , \, A}  \, .
  \label{exesr1}
\end{equation}
\end{theorem}  

\begin{corollary} \label{exerrest}
For $j \in \mathbb{N}$ and $0 \le \mu \le j$,  if $D^{j}f \in L^{2}_{\rho^{(\beta + j - 1 , \alpha - \beta + j - 1)}}(0,1)$, then
there exists $C > 0$ (independent of $N$ and $\alpha$) such that
\begin{equation}
 \|(q \, - \, q_{N} )/\rho^{(\alpha-\beta , \beta)}\|_{\mu , \, \rho^{(\alpha - \beta \, , \, \beta)} , \, A} \ 
 \le \ C \, N^{- (\alpha - 1)} \,
  ( N \, (N + \alpha - 2) )^{\frac{\mu - j}{2}} \, | f |_{j , \, \rho^{(\beta - 1 \, , \, \alpha - \beta - 1)} , \, A}  \, .
  \label{hesr1}
\end{equation}
\end{corollary}
\textbf{Proof}:
Noting that $f_{N}(x) \ = \ P_{N , \, \beta - 1 \, , \, \alpha - \beta - 1}f(x)$, from \eqref{thy1}, taking the limit as
$M \rightarrow \infty$, we have
\begin{align}
\| D^{\mu}\big( (q \, - \, q_{N}) / \rho^{(\alpha-\beta , \beta)} \big)\|_{\rho^{(\alpha - \beta + \mu \, , \, \beta + \mu)}} 
&\le C \, N^{- (\alpha - 1)} \, \| D^{\mu} (f \, - \, f_{N}) \|_{\rho^{(\beta + \mu - 1 \, , \, \alpha - \beta + \mu - 1)}}    \nonumber \\
&\le C \, N^{- (\alpha - 1)} \, \| f \, - \, f_{N} \|_{\mu , \, \rho^{(\beta  - 1 \, , \, \alpha - \beta  - 1)} , \, A}    \nonumber \\
&\le C \, N^{- (\alpha - 1)} \, ( N \, (N + \alpha - 2) )^{\frac{\mu - j}{2}} \, | f |_{j , \, \rho^{(\beta - 1 \, , \, \alpha - \beta - 1)}  , \, A} \, ,
\end{align}
where, in the last step we have used \eqref{exesr1}). 
\mbox{ } \hfill \qed 

\setcounter{equation}{0}
\setcounter{figure}{0}
\setcounter{table}{0}
\setcounter{theorem}{0}
\setcounter{lemma}{0}
\setcounter{corollary}{0}
\section{Convergence of $( u(x) - u_{N}(x) )$}
\label{ssecEc3}
From \eqref{FDE:e2}, $u(x)$ is given by
\[
 u(x) \ = \ - \, \int_{0}^{x} \frac{\tilde{q}(s)}{K(s)} \, ds \ = \ - \, c_{-2} \int_{0}^{x} \frac{k(s)}{K(s)} \, ds \ - \ 
 \int_{0}^{x} \frac{q(s)}{K(s)} \, ds \, . 
 \]
Hence,
\begin{align}
| u(x) - u_{N}(x) | &\le \ | c_{-2} - c_{-2 , N} | \, \left|  \int_{0}^{x} \frac{k(s)}{K(s)} \, ds \right| \ + \ 
 \left| \int_{0}^{x} \frac{q(s) - q_{N}(s)}{K(s)} \, ds \right|   \nonumber \\
&\le \
   \ | c_{-2} - c_{-2 , N} | \, \left|  \int_{0}^{1} \frac{k(s)}{K(s)} \, ds \right|   \nonumber \\
   & \quad \quad \ + \ 
    \left( \int_{0}^{1} (1 - s)^{-(\alpha - \beta)} \, s^{- \beta} \, ( q(s) - q_{N}(s) )^{2} \, ds \right)^{1/2} \, 
     \left( \int_{0}^{1} \frac{(1 - s)^{\alpha - \beta} \, s^{\beta}}{K^{2}(s)} \, ds \right)^{1/2}   \nonumber \\
&\le \
   C_{11} \, C_{12} \, \| q - q_{N} \|_{\rho^{(-(\alpha - \beta) , - \beta)}}  \, C_{11}^{-1} \ + \ 
     C_{12} \, \| q - q_{N} \|_{\rho^{(-(\alpha - \beta) , - \beta)}}    \nonumber \\
&=  \ 2 \, C_{12} \, \| q - q_{N} \|_{\rho^{(-(\alpha - \beta) , - \beta)}}  \, .  \label{cvgu}
\end{align}     
   
The above analysis is very coarse and most likely does not give the best error estimate for $( u(x) - u_{N}(x) )$.
In the next section we obtain a better error estimate for the special case of $K(x) = $ constant.

\subsection{Convergence of $( u(x) - u_{N}(x) )$ -- Special case $K(x) = $ constant}
\label{ssecEc4}
In this section we investigate the convergence of $u_{N}(x)$ to $u(x)$ for the special case when $K(x) = $ constant.
\begin{corollary} \label{erruuNS}
For $j \in \mathbb{N}$ and $0 \le \mu \le j$,  if $D^{j}f \in L^{2}_{(1 - x)^{\beta + j - 1} \, x^{\alpha - \beta + j - 1}}(0,1)$, then
there exists $C > 0$ (independent of $N$ and $\alpha$) such that
\begin{equation}
 \|  u \, - \, u_{N} \|_{\rho^{(-(\alpha - \beta + 1) \, , \, -(\beta + 1))}} \ \leq \ C \, 
 N^{-1} \, (N + 1)^{-(\alpha - 1)} \, (N \, (N + \alpha - 1))^{-\frac{j}{2}} \, | f |_{j,\rho^{(\beta-1,\alpha-\beta-1)},A}.
 \label{erstu12}
\end{equation}
\end{corollary}
\textbf{Proof}:
Recall that $\tilde{q}(x) \ = \ c_{-2} k(x) \ + \ q(x)$, where $c_{-2}$ is determined by \eqref{FDE:e3}. When $K(x) = $ constant,
\begin{equation}
 c_{-2} \ = \ - \frac{ \int_{0}^{1} \frac{q(s)}{K(s)} \, ds}{  \int_{0}^{1} \frac{k(s)}{K(s)} \, ds}
  \ = \ - \frac{ \int_{0}^{1} q(s) \, ds}{ \int_{0}^{1} k(s) \, ds} \, .
 \label{calc2}
\end{equation}

From \eqref{eqC2} it follows that
\begin{equation}
 \int_{0}^{x} \rho^{(\alpha - \beta \, \beta)}(s) \, G_{n}^{(\alpha - \beta \, \beta)}(s) \, ds 
 \ = \ \frac{-1}{n}  \rho^{(\alpha - \beta + 1 \, \beta + 1)}(x) \, G_{n - 1}^{(\alpha - \beta + 1 \, \beta + 1)}(x) \, , \quad n \ge 1\, .
 \label{eqpp2}
\end{equation}

Now,
\begin{align}
u_{M}(x) \, - \, u_{N}(x) &= \int_{0}^{x} \rho^{(\alpha - \beta \, , \, \beta)}(s) \, 
 \sum_{i = N}^{M-1}  c_{i} \, G_{i}^{(\alpha - \beta \, , \, \beta)}(s) \, ds  \nonumber \\
 &= \ - \rho^{(\alpha - \beta + 1 \, , \, \beta + 1)}(x) \, 
 \sum_{i = N}^{M-1}  \frac{1}{i} c_{i} \, G_{i-1}^{(\alpha - \beta + 1\, , \, \beta + 1)}(x) \, , \quad \mbox{(using \eqref{eqpp2})} \, .
 \nonumber 
\end{align} 

Hence,
\begin{align}
& \| u_{M} \, - \, u_{N} \|_{\rho^{(-(\alpha - \beta + 1) \, , \, -(\beta + 1))}}^{2}   \nonumber \\
 & \quad = \ 
  \bigg( \rho^{(\alpha - \beta + 1 \, , \, \beta + 1)}(x) \, 
 \sum_{i = N-1}^{M-2}  \frac{c_{i+1} \, G_{i}^{(\alpha - \beta + 1 \, , \, \beta + 1)}(x)}{i+1}  \, , \,
  \sum_{i = N-1}^{M-2}  \frac{c_{i+1} \, G_{i}^{(\alpha - \beta + 1 \, , \, \beta + 1)}(x)}{i+1}  \bigg)  \nonumber \\
& \quad = \
   \sum_{i = N-1}^{M-2}  \frac{1}{(i+1)^{2}} \frac{f_{i+2}^{2}}{\lambda_{i+1}^{2} \,  |\| G_{i+2}^{(\beta - 1 \, , \, \alpha - \beta - 1)} |\|^{4}} \,
   |\| G_{i}^{(\alpha - \beta + 1 \, , \, \beta + 1)} |\|^{2}  \nonumber \\
& \quad = \
   \sum_{i = N+1}^{M}  \frac{1}{(i-1)^{2}} \frac{f_{i}^{2}}{\lambda_{i-1}^{2} \,  |\| G_{i}^{(\beta - 1 \, , \, \alpha - \beta - 1)} |\|^{4}} \,
   |\| G_{i-2}^{(\alpha - \beta + 1 \, , \, \beta + 1)} |\|^{2}  \label{ytr45}
\end{align}   

Similar to Lemma \ref{lmageq}, we have
\begin{align}
   \frac{  |\| G_{i-2}^{(\alpha - \beta + 1 \, , \, \beta + 1)} |\|^{2}}{|\| G_{i}^{(\beta - 1 \, , \, \alpha - \beta - 1)} |\|^{2}} \ &= \frac{1}{2i+\alpha-1}  
   \, \frac{\Gamma(i + \alpha - \beta) \, \Gamma(i + \beta)}{\Gamma(i - 1) \Gamma(i + \alpha + 1)} \, \frac{2i + \alpha -1}{1} \, 
   \frac{\Gamma(i + 1) \, \Gamma(i + \alpha-1)}{\Gamma(i + \beta) \, \Gamma(i + \alpha - \beta)}   \nonumber \\
   &=\frac{i (i-1)}{(i + \alpha)(i + \alpha - 1)} \leq 1 \, .  \label{yytt22}
\end{align}
Using \eqref{ytr45} and \eqref{yytt22}, together with (\ref{nwbde}) we then obtain 
\begin{align*}
&\| u_{M} \, - \, u_{N} \|_{\rho^{(-(\alpha - \beta + 1) \, , \, -(\beta + 1))}}^{2} 
\le   \sum_{i = N+1}^{M}  \frac{1}{(i-1)^{2}} \frac{f_{i}^{2}}{\lambda_{i-1}^{2} \,  |\| G_{i}^{(\beta - 1 \, , \, \alpha - \beta - 1)} |\|^{2}} \\
& \quad  \leq \ \frac{1}{N^{2} \, \lambda_{N}^{2}} \left(\rho^{(\beta - 1 \, , \, \alpha - \beta - 1)} \, \sum_{i = N+1}^{M} 
 \frac{f_{i}}{ |\| G_{i}^{(\beta - 1 \, , \, \alpha - \beta - 1)} |\|^{2}} \, G_{i}^{(\beta - 1 \, , \, \alpha - \beta - 1)}(x) \, , \right. \\
& \quad \quad \quad \quad \quad \quad  \quad \quad \quad  \quad \quad \quad \quad \quad \quad  \quad \quad \quad 
\left. \sum_{i = N+1}^{M}  \frac{f_{i}}{ |\| G_{i}^{(\beta - 1 \, , \, \alpha - \beta - 1)} |\|^{2}} \,
 G_{i}^{(\beta - 1 \, , \, \alpha - \beta - 1)}(x) \right)\\
& \quad  \leq \ \frac{C}{N^{2} \, (N+1)^{2(\alpha - 1)}} \, \| f_{M} \, - \, f_{N} \|_{\rho^{(\beta - 1,\alpha - \beta -1) }}^{2}.
\end{align*}

Then, similar to the proof of Corollary \ref{exerrest} we get
\begin{equation}
 \|  u \, - \, u_{N} \|_{\rho^{(-(\alpha - \beta + 1) \, , \, -(\beta + 1))}} \ \leq \ C \, 
 N^{-1} \, (N + 1)^{-(\alpha - 1)} \, (N \, (N + \alpha - 1))^{-\frac{j}{2}} \, | f |_{j,\rho^{(\beta-1,\alpha-\beta-1)},A}.
\end{equation} 
\mbox{ } \hfill \qed

%%%%
%%%%
%%%%
\setcounter{equation}{0}
\setcounter{figure}{0}
\setcounter{table}{0}
\setcounter{theorem}{0}
\setcounter{lemma}{0}
\setcounter{corollary}{0}
\section{Numerical experiments}
\label{secNum}
In this section we present two numerical examples to demonstrate our approximation scheme, and to compare
the experimental rate of convergence of the approximation with the theoretically predicated rate. Within 
Example 1 we consider three numerical experiments corresponding to different values of $\alpha$ and $r$. 
For this example we choose $K(x) = 1$ which permits us to compare the theoretically predicted rate
of convergence of $u_{N}$ to $u$ in the $L^2_{\rho^{(-(\alpha-\beta+1),-(\beta+1))}}$ norm with 
its experimental rate.

In order to determine the theoretical rate of convergence for $\|q-q_N\|_{L^2_{\rho^{(-(\alpha-\beta),-\beta)}}}$
and \linebreak[4]
$\|u-u_N\|_{L^2_{\rho^{(-(\alpha-\beta+1),-(\beta+1))}}}$ from \eqref{hesr1} and \eqref{erstu12}, respectively,
we need to determine the largest value for $j$ such that $f(x) \in H^{j}_{\rho^{(\beta - 1 \, , \, \alpha - \beta - 1)} \,  , \, A}(0,1)$,
i.e, the largest $j$ such that $\| D^{j} f \|^{2}_{\rho^{(\beta + j - 1 \, , \, \alpha + j - \beta - 1)}} \, < \, \infty$.
The most singular terms for $f(x)$ in Example 1 are $x^{2 - \alpha}$ and $(1 - x)^{2 - \alpha}$. We focus our attention
on $x^{2 - \alpha}$.

Note that $D^{j} x^{2 - \alpha} \sim x^{2 - \alpha - j}$. Then 
\begin{align*}
\| D^{j} x^{2 - \alpha} \|^{2}_{\rho^{(\beta + j - 1 \, , \, \alpha + j - \beta - 1)}} 
&\sim \ \int_{0}^{1} x^{\alpha + j - \beta - 1} \, \left( x^{2 - \alpha - j} \right)^{2} \, dx  
\ = \ \int_{0}^{1} x^{3 - \alpha -  \beta - j } \, dx \ < \ \infty  \\
\Rightarrow \ \ -1 &<  \ 3 - \alpha -  \beta - j  \\
\Rightarrow \ \  j &<  \ 4 - \alpha -  \beta \, .  
\end{align*}

Then, for experiment 1 in Example 1 ($\alpha = 1.60$, $\beta = 0.85$) 
$f(x) \in H^{j}_{\rho^{(\beta - 1 \, , \, \alpha - \beta - 1)} \,  , \, A}(0,1)$ for $j < 1.55$, which leads to 
theoretical asymptotic rates of $\|q-q_N\|_{L^2_{\rho^{(-(\alpha-\beta),-\beta)}}} \sim  \, N^{-(\alpha - 1 + j)} 
\, =  \, N^{- 2.15}$
and $\|u-u_N\|_{L^2_{\rho^{(-(\alpha-\beta+1),-(\beta+1))}}} \sim  \, N^{-(\alpha + j)} 
\, = \,  N^{- 3.15}$.

Assuming that $\| \xi - \xi_{N} \|_{L_{\rho}} \sim  \, N^{- \kappa}$, the experimental convergence rate
is calculated using
\[
  \kappa \approx \frac{\log ( \| \xi - \xi_{N_{1}} \|_{L_{\rho}}  /  \| \xi - \xi_{N_{2}} \|_{L_{\rho}} )}{ \log (N_{2} / N_{1})} \, .
\]

\textbf{Example 1}. Let $K(x)=1$, $\beta$ be determined by \eqref{defbeta}, and
\begin{align*}
 f(x)&=\frac{6r}{\Gamma(2-\alpha)\delta}((2\alpha-8)x^{3-\alpha}+(\alpha-3)(\alpha-4)x^{2-\alpha})\\[0.025in]
&+\frac{6(1-r)}{\Gamma(2-\alpha)\delta}(-(2\alpha-8)(1-x)^{3-\alpha}-(\alpha-3)(\alpha-4)(1-x)^{2-\alpha}),
\end{align*}
where $\delta:=\alpha^3-9\alpha^2+26\alpha-24$. Then the solution $u(x)$, and the related $q(x)$, are given by
\begin{align*}
u(x) \, = \, 3x^{2} \, - \, 2x^{3} \, - \, 
\frac{x^{\beta} \, {}_2 F_1(-\alpha + \beta + 1 \, , \, \beta \, ; \, \beta + 1 \, , \, x)}%
{{}_2 F_1(-\alpha + \beta + 1 \, , \, \beta \, ; \, \beta + 1 \, , \, 1)},~~q(x) \, = \, -6x + 6x^{2},
\end{align*}
where ${}_2 F_1(a,b;c,x)$ donate the Gauss three-parameter hypergeometric function defined by an integral and series as follows:
 \begin{equation*}
  \begin{array}{rcl}
 {}_2 F_1(a,b;c,x)& = &\displaystyle \frac{\Gamma (c)}{\Gamma (b)\Gamma (c-b)} \int_{0}^{1}z^{b-1}(1-z)^{c-b-1}(1-zx)^{-a}dz \\
 & = & \displaystyle \sum\limits_{n=0}^{\infty} \frac{(a)_n(b)_nx^n}{(c)_nn!}, 
 \end{array}
 \end{equation*}
 with convergence only if $Re(c)>Re(b)>0$ and $(s)_n$ is the rising Pochhammer symbol defined by 
 $(s)_n=\Gamma (s+n)/\Gamma (s)$. 

\begin{figure}[H]
	\setlength{\abovecaptionskip}{0pt}
	\centering
	\includegraphics[width=2.8in,height=2.8in]{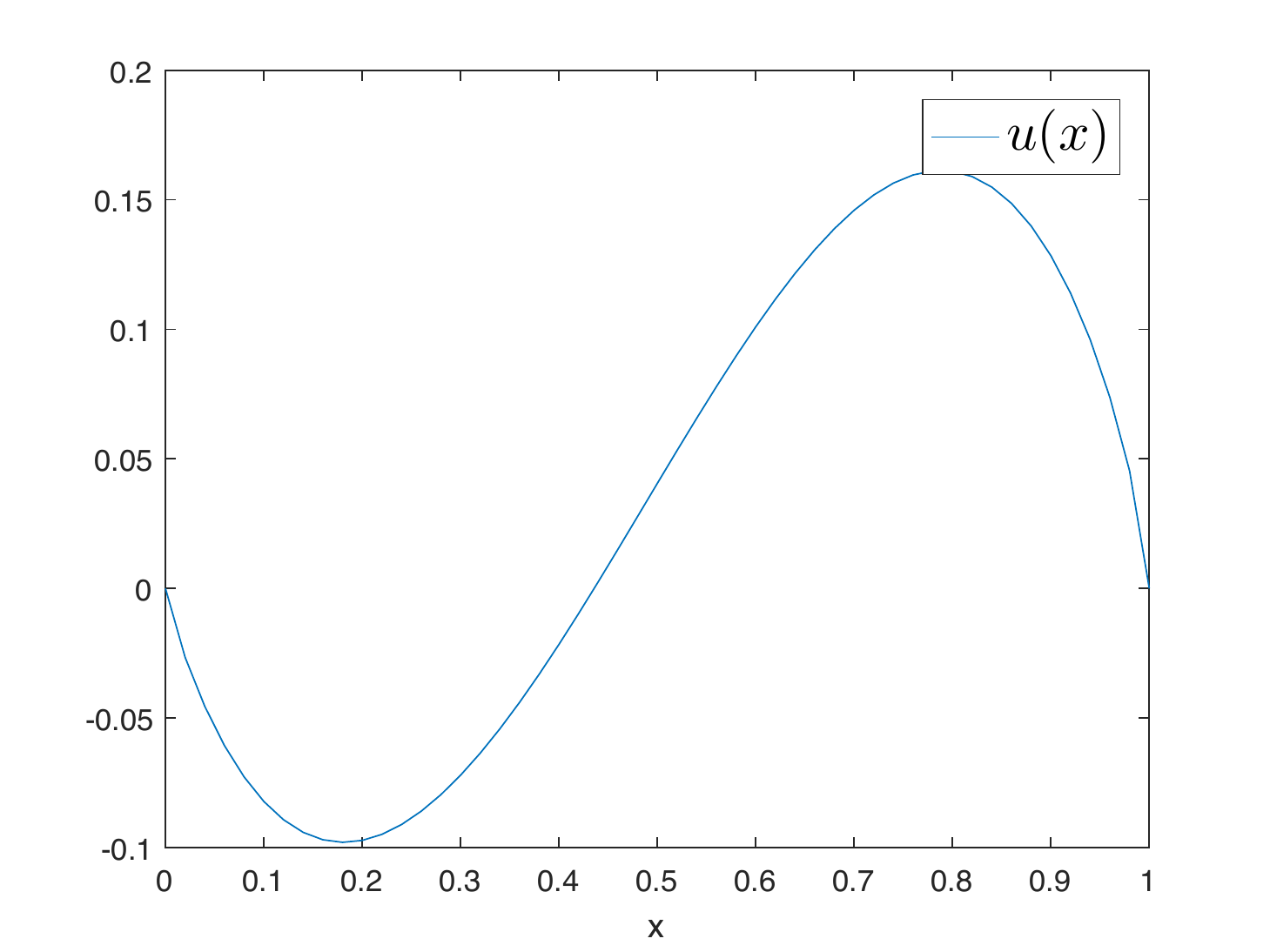}
	\includegraphics[width=2.8in,height=2.8in]{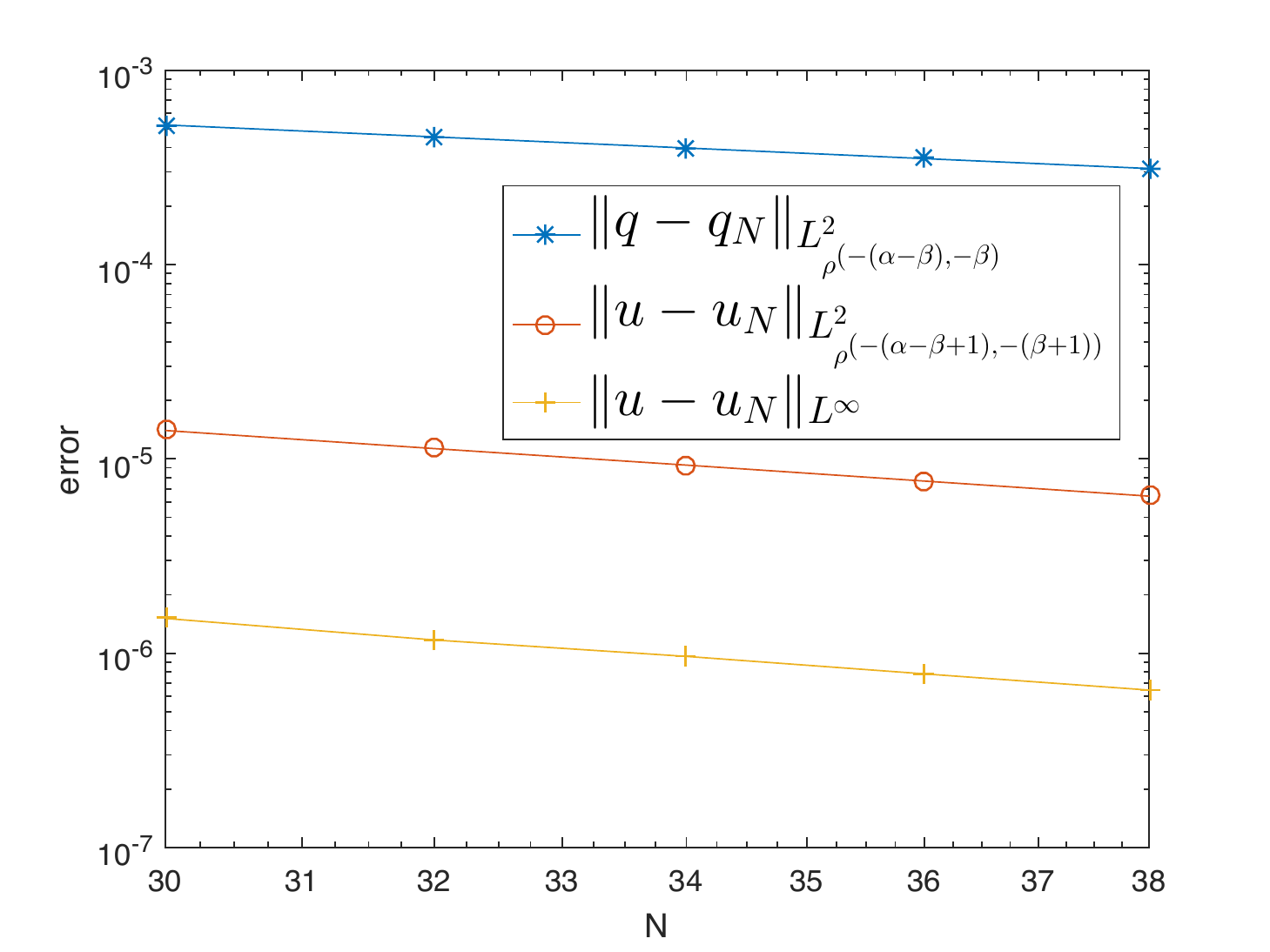}
	\caption{The plot of solution $u(x)$ (left), and (right) the  $log-log$ plot of the errors for experiment 1 of Example 1.}
	\label{logloggraph}
\end{figure}

A plot of the solution $u(x)$, corresponding to $\alpha = 1.60$, $r = 0.39$ and $\beta = 0.85$, and a plot of 
the errors for this numerical experiment are presented in Figure \ref{logloggraph}.

The experimental convergence rate $\kappa$ of the error in different norms for Example 1 are shown 
in Table \ref{1,1}, \ref{1,2} and \ref{1,3}.

 \begin{table}[H]
 	\setlength{\abovecaptionskip}{0pt}
 	\centering
 	\caption{Example 1 with $\alpha = 1.60$, $r = 0.39$ and $\beta = 0.85$.}	\label{1,1}
 	\vspace{0.5em}	
 	\begin{tabular}{ccccccc}
 		\hline
  $N$ & $\|q-q_N\|_{L^2_{\rho^{(-(\alpha-\beta),-\beta)}}}$ & $\kappa$ & $\|u-u_N\|_{L^2_{\rho^{(-(\alpha-\beta+1),-(\beta+1))}}}$ 
  & $\kappa$ & $\|u-u_N\|_{L^{\infty}}$ & $\kappa$\\
 		\cline{1-7}
30&	5.23E-04&		&    1.40E-05&		&    1.51E-06&	\\
32& 4.54E-04&	2.18 &	1.13E-05 &  3.26& 	1.17E-06&	3.90 \\
34&	3.98E-04&	2.18 &	9.29E-06&	3.28 &	9.63E-07&	3.27 \\
36&	3.51E-04&	2.18 &	7.69E-06&	3.30 &	7.83E-07&	3.62 \\
38&	3.12E-04&	2.18 &	6.43E-06&	3.33 &	6.46E-07&	3.54 \\
 	\hline
 	Pred.&   &2.15&   &3.15&   &  2.15\\
 		\hline	
 	\end{tabular}
 \end{table}
 
\begin{table}[H]
 	\setlength{\abovecaptionskip}{0pt}
 	\centering
 	\caption{Example 1 with $\alpha = 1.40$, $r = 0.50$ and $\beta = 0.70$.}	\label{1,2}
 	\vspace{0.5em}	
 	\begin{tabular}{ccccccc}
 		\hline
  $N$ & $\|q-q_N\|_{L^2_{\rho^{(-(\alpha-\beta),-\beta)}}}$ & $\kappa$ & $\|u-u_N\|_{L^2_{\rho^{(-(\alpha-\beta+1),-(\beta+1))}}}$ 
  & $\kappa$ & $\|u-u_N\|_{L^{\infty}}$ & $\kappa$\\
 		\cline{1-7}
30&  5.10E-04&		&1.37E-05&		&1.59E-06&	\\
32&  4.40E-04&	2.29 &	1.10E-05&	3.39& 	1.22E-06&	4.16\\ 
34&  3.83E-04&	2.29 &	8.94E-06&	3.42 &	1.02E-06&	2.87 \\
36&  3.36E-04&	2.29 &	7.34E-06&	3.44 &	8.28E-07&	3.72 \\
38&  2.97E-04&	2.29 &	6.09E-06&	3.47 &	6.72E-07&	3.86 \\
 	\hline
 	Pred.&   &2.30&   &3.30&   & 2.30 \\
 		\hline	
 	\end{tabular}
 \end{table}
 
 \begin{table}[H]
 	\setlength{\abovecaptionskip}{0pt}
 	\centering
 	\caption{Example 1 with $\alpha = 1.80$, $r = 0.50$ and $\beta = 0.90$.}	\label{1,3}
 	\vspace{0.5em}	
 	\begin{tabular}{ccccccc}
 		\hline
  $N$ & $\|q-q_N\|_{L^2_{\rho^{(-(\alpha-\beta),-\beta)}}}$ & $\kappa$ & $\|u-u_N\|_{L^2_{\rho^{(-(\alpha-\beta+1),-(\beta+1))}}}$ 
  & $\kappa$ & $\|u-u_N\|_{L^{\infty}}$ & $\kappa$\\
 		\cline{1-7}
30&	4.21E-04&		&1.11E-05&		&1.08E-06&	\\
32&	3.69E-04&	2.07 &	9.07E-06&	3.16& 	8.40E-07&	3.85\\ 
34& 3.25E-04&	2.07 &	7.48E-06&	3.18 &	7.08E-07&	2.82 \\
36&	2.89E-04&	2.08 &	6.23E-06&	3.21 &	5.76E-07&	3.60 \\
38&	2.58E-04&	2.08 &	5.23E-06&	3.24 &	4.64E-07&	4.03 \\
 	\hline
 	Pred.&  &2.10&  &3.10&  & 2.10 \\
 		\hline	
 	\end{tabular}
 \end{table}
 
The experimental convergence rates for $\|q-q_N\|_{L^2_{\rho^{(-(\alpha-\beta),-\beta)}}}$ and
$\|u-u_N\|_{L^2_{\rho^{(-(\alpha-\beta+1),-(\beta+1))}}}$ are in strong agreement  with the theoretically predicted rates. Not 
surprisingly, the theoretically predicted rate for $\|u-u_N\|_{L^{\infty}}$ in \eqref{cvgu} appears to be suboptimal.

\textbf{Example 2}. With this example we investigate the numerical approximation for the interesting case
of a non constant $K(x)$. Let $K(x)=1+x^2$ and
\begin{align*}
 f(x)&=r\left(-480\frac{x^{6-\alpha}}{\Gamma(7-\alpha)}+144\frac{x^{5-\alpha}}{\Gamma(6-\alpha)}-36\frac{x^{4-\alpha}}{\Gamma(5-\alpha)}+12\frac{x^{3-\alpha}}{\Gamma(4-\alpha)}-2\frac{x^{2-\alpha}}{\Gamma(3-\alpha)}\right)\\
 &-(1-r)\left(480\frac{(1-x)^{6-\alpha}}{\Gamma(7-\alpha)}-366\frac{(1-x)^{5-\alpha}}{\Gamma(6-\alpha)}+132\frac{(1-x)^{4-\alpha}}{\Gamma(5-\alpha)}-32\frac{(1-x)^{3-\alpha}}{\Gamma(4-\alpha)}+4\frac{(1-x)^{2-\alpha}}{\Gamma(3-\alpha)}\right).
\end{align*}
Then the solution $u(x)$, and the related $q(x)$, are 
\begin{align*}
u(x) \, = \, x^{2} (1-x)^{2},~~q(x) \, = \, - 2 (1 + x^{2}) x (1 - x) (1 - 2x).
\end{align*}
The convergence rate $\kappa$ of the error in different norms for Example 2 are shown in Table \ref{2,1}, \ref{2,2} and \ref{2,3}.
The numbers given for the predicted rate of convergence of 
$\|u-u_N\|_{L^2_{\rho^{(-(\alpha-\beta+1),-(\beta+1))}}}$, denoted with an $\mbox{}^{*}$, are from \eqref{erstu12}, which does not apply in this setting as
$K(x) \ne constant$.

\begin{table}[H]
 	\setlength{\abovecaptionskip}{0pt}
 	\centering
 	\caption{Example 2 with $\alpha = 1.60$, $r = 0.39$ and $\beta = 0.85$.}	\label{2,1}
 	\vspace{0.5em}	
 	\begin{tabular}{ccccccc}
 		\hline
 		$N$&$\|q-q_N\|_{L^2_{\rho^{(-(\alpha-\beta),-\beta)}}}$ &$\kappa$&$\|u-u_N\|_{L^2_{\rho^{(-(\alpha-\beta+1),-(\beta+1))}}}$ &$\kappa$&$\|u-u_N\|_{L^{\infty}}$&$\kappa$\\
 		\cline{1-7}
30&	3.01E-04&		&5.57E-06	&	&7.21E-07&	\\
32&	2.59E-04&	2.28& 	4.50E-06&	3.31& 	5.54E-07&	4.07\\ 
34&	2.26E-04&	2.28 &	3.68E-06&	3.30 &	4.54E-07&	3.28 \\
36&	1.98E-04&	2.28 &	3.05E-06&	3.28 &	3.63E-07&	3.95 \\
38&	1.75E-04&	2.27 &	2.56E-06&	3.27 &	2.92E-07&	4.01 \\
 	\hline
 	Pred.&    &2.15&   &$3.15^{*}$&   &  2.15  \\
 		\hline	
 	\end{tabular}
 \end{table}
 
\begin{table}[H]
 	\setlength{\abovecaptionskip}{0pt}
 	\centering
 	\caption{Example 2 with $\alpha = 1.40$, $r = 0.50$ and $\beta = 0.70$.}	\label{2,2}
 	\vspace{0.5em}	
 	\begin{tabular}{ccccccc}
 		\hline
 		$N$&$\|q-q_N\|_{L^2_{\rho^{(-(\alpha-\beta),-\beta)}}}$ &$\kappa$&$\|u-u_N\|_{L^2_{\rho^{(-(\alpha-\beta+1),-(\beta+1))}}}$ &$\kappa$&$\|u-u_N\|_{L^{\infty}}$&$\kappa$\\
 		\cline{1-7}
30&	2.90E-04&		&5.49E-06&		&7.73E-07&	\\
32&	2.49E-04&	2.37 &	4.40E-06&	3.42& 	6.08E-07&	3.72\\ 
34&	2.16E-04&	2.36 &	3.58E-06&	3.41 &	4.79E-07&	3.92 \\
36&	1.89E-04&	2.36 &	2.95E-06&	3.40 &	3.82E-07&	3.97 \\
38&	1.66E-04&	2.35 &	2.45E-06&	3.39 &	3.10E-07&	3.83 \\
 	\hline
 	Pred.&  &2.30&   &$3.30^{*}$&    & 2.30 \\
 		\hline	
 	\end{tabular}
 \end{table}
 
 \begin{table}[H]
 	\setlength{\abovecaptionskip}{0pt}
 	\centering
 	\caption{Example 2 with $\alpha = 1.80$, $r = 50$ and $\beta = 0.90$.}	\label{2,3}
 	\vspace{0.5em}	
 	\begin{tabular}{ccccccc}
 		\hline
 		$N$&$\|q-q_N\|_{L^2_{\rho^{(-(\alpha-\beta),-\beta)}}}$ &$\kappa$&$\|u-u_N\|_{L^2_{\rho^{(-(\alpha-\beta+1),-(\beta+1))}}}$ &$\kappa$&$\|u-u_N\|_{L^{\infty}}$&$\kappa$\\
 		\cline{1-7}
30&	2.38E-04&		&4.40E-06&		&5.22E-07&	\\
32&	2.07E-04&	2.14 &	3.58E-06&	3.19& 	4.10E-07&	3.73\\ 
34&	1.82E-04&	2.14 &	2.95E-06&	3.18 &	3.32E-07&	3.49 \\
36&	1.61E-04&	2.14 &	2.46E-06&	3.17 &	2.68E-07&	3.77 \\
38&	1.43E-04&	2.14 &	2.08E-06&	3.16 &	2.16E-07&	4.00 \\
 	\hline
 	Pred.&   &2.10&   &$3.10^{*}$&   & 2.10\\
 		\hline	
 	\end{tabular}
 \end{table}
 
The experimental convergence rate for $\|q-q_N\|_{L^2_{\rho^{(-(\alpha-\beta),-\beta)}}}$ is in good agreement with
that predicted theoretically. Though estimate \eqref{erstu12} does not apply for $K(x) \ne constant$, nonetheless the
predicted value using \eqref{erstu12} is in good agreement with the experimental convergence rate. We again note that
the error estimate obtained in \eqref{cvgu} appears to be suboptimal.
 
 \section*{Acknowledgements}
 
 This work was funded by the OSD/ARO MURI Grant W911NF-15-1-0562 and by the National Science Foundation under Grant DMS-1620194.

\end{document}